\documentclass[a4paper,12pt,oneside]{article}
\usepackage{amsmath,amssymb,amsfonts,graphicx}

\newtheorem{theorem}{Theorem}[section]

\newtheorem{corollary}[theorem]{Corollary}
\newtheorem{lemma}[theorem]{Lemma}

\begin{document}

\def\Z{\mathbb{ Z}}
\def\C{\mathbb{ C}}
\def\F{\mathbb{ F}}
\def\R{\mathbb{ R}}
\def\E{\mathbb{ E}}
\def\S{\mathbb{ S}}
\def\P{\mathbb{ P}}
\def\ra{\rightarrow}
\def\l{\ell}

\title{The Biham-Middleton-Levine traffic model\\ for a single junction}
\author{I. Benjamini \quad O. Gurel-Gurevich \quad R. Izkovsky}
\maketitle
\begin{abstract}
In the Biham-Middleton-Levine traffic model cars are placed in
some density $p$ on a two dimensional torus, and move according to
a (simple) set of predefined rules. Computer simulations show this
system exhibits many interesting phenomena: for low densities the
system self organizes such that cars flow freely while for
densities higher than some critical density the system gets stuck
in an endless traffic jam. However, apart from the simulation
results very few properties of the system were proven rigorously
to date. We introduce a simplified version of this model in which
cars are placed in a single row and column (a junction) and show
that similar phenomena of self-organization  of the system and
phase transition still occur.
\end{abstract}

\section{The BML traffic model}
\paragraph{}
The Biham-Middleton-Levine (BML) traffic models was first
introduced in \cite{bml} published 1992. The model involves two
types of cars: "red" and "blue". Initially the cars are placed in
random with a given density $p$ on the $N\times N$ torus. The
system dynamics are as follows: at each turn, first all the red
cars try to move simultaneously a single step to the right in the
torus. Afterwards all blue cars try to move a single step upwards.
A car succeeds in moving as long as the relevant space
above/beside it (according to whether it is blue/red) is vacant.

\begin{figure}
\centering
\includegraphics[scale=0.5]{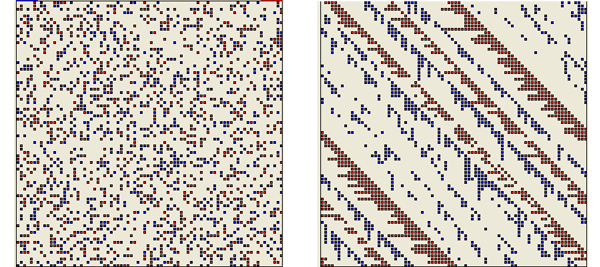}
\caption{Self organization for $p=0.3$: the initial configuration
in the left organizes after 520 steps to the "diagonal"
configuration on the right which flows freely} \label{cars1}
\end{figure}

\begin{figure}
\centering
\includegraphics[scale=0.5]{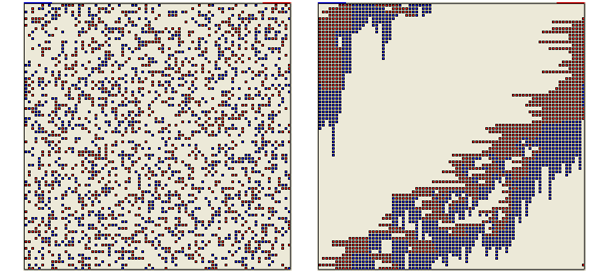}
\caption{Traffic jam for $p=0.35$: the initial configuration on
the left gets to the stuck configuration on the right after 400
steps, in which no car can ever move} \label{cars2}
\end{figure}

\paragraph{}
The basic properties of this model are described in \cite{bml} and
some recent more subtle observations due to larger simulations are
described in \cite{dsouza}. The main and most interesting property
of the system originally observed in simulations is a \emph{phase
transition}: for some critical density $p_c$ one observes, that
while filling the torus randomly with cars in density $p < p_c$
the system self organizes such that after some time all cars flow
freely with no car ever running into another car (see figure
\ref{cars1}), by slightly changing the density to some $p > p_c$
not only does system not reach a free flow, but it will eventually
get stuck in a configuration in which no car can ever move (see
figure \ref{cars2}).
\paragraph{}
Very little of the above behaviour is rigorously proven for the BML
model. The main rigorous result is that of Angel, Holroyd and Martin
\cite{omer}, showing that for some fixed density $p<1$, very close
to 1, the probability of the system getting stuck tends to 1 as
$N\ra \infty$. In \cite{itai} one can find a study of the very low
density regime (when $p = O(\frac{1}{N})$).

\paragraph{}
First, we introduce a slight variant of the original BML model by
allowing a car (say, red) to move not only if there is a vacant
place right next to it but also if there is a red car next to it
that moves. Thus, for sequence of red cars placed in a row with a
single vacant place to its right - all cars will move together (as
oppose to only the rightmost car in the sequence for the original
BML model). Not only does this new variant exhibits the same
phenomena of self-organization and phase transition, they even seem
to appear more quickly (i.e. it takes less time for the system to
reach a stable state). Actually the demonstrated simulations
(figures \ref{cars1}, \ref{cars2}) were performed using the variant
model. Note that the results of \cite{omer} appear to apply equally
well to the variant model.

\paragraph{}
In the following, we will analyze a simplified version of the BML
model: BML on a single junction. Meaning, we place red cars in some
density $p$ on a single row of the torus, and blue cars are placed
in density $p$ on a single column. We will show that for all $p$ the
system reaches \emph{optimal speed}\footnote{The speed of the system
is the asymptotic average rate in which a car moves - i.e. number of
actual moves per turn.}, depending only on $p$. For $p < 0.5$ we
will show the system reaches speed 1, while for $p>0.5$ the speed
cannot be 1, but the system will reach the same speed, regardless of
the initial configuration. Moreover, at $p=0.5$ the system's
behaviour undergoes a phase transition: we will prove that while for
$p < 0.5$ the stable configuration will have linearly many sequences
of cars, for $p > 0.5$ we will have only $O(1)$ different sequences
after some time. We will also examine what happens at a small window
around $p=0.5$.

\paragraph{}
Note that in the variant BML model (and unlike the original BML
model) car sequences are never split. Therefore, the simplified
version of the variant BML model can be viewed as some kind of
1-dimensional coalescent process.

\paragraph{}
Much of the proofs below rely on the fact that we model BML on a
symmetric torus. Indeed, the time-normalization of section
\ref{normalization} would take an entirely different form if the
height and width were not equal. We suspect that the model would
exhibit similar properties if the height and width had a large
common denominator, e.g. if they had a fixed proportion. More
importantly, we believe that for height and width with low common
denominator (e.g. relatively prime), we would see a clearly
different behaviour. As a simple example, note that in the
relatively prime case, a speed of precisely 1 cannot be attained, no
matter how low is $p$, in contrast with corollary \ref{speed1cor}.
This dependence on the arithmetic properties of the dimensions is
also apparent in \cite{dsouza}.

\section{The Junction model}

\subsection{Basic model}

\paragraph{}
We start with the exact definition of our simplified model. On a
cross shape, containing a single horizontal segment and a single
vertical segment, both of identical length $N$, red cars are placed
in exactly $pN$ randomly (and uniformly) chosen locations along the
row, and blue cars are similarly placed in $pN$ locations along the
column. $p$ will be called the density of the configuration. For
simplicity and consistency with the general BML model, we refer to
the cars placed on the horizontal segment as red and those on the
vertical segment as blue. For simplicity we may assume that the
\emph{junction}, i.e. the single location in the intersection of the
segments, is left unoccupied. The segments are cyclic - e.g. a red
car located at the rightmost position ($N-1$) that moves one step to
the right re-emerges at the leftmost position (0).

\paragraph{}
At each turn, all the red cars move one step to the right, except
that the red car that is just left of the junction will not move if
a blue car is in the junction (i.e. blocking it), in which case also
all the red cars immediately following it will stay still.
Afterwards, the blue cars move similarly, with the red and blue
roles switched. As in the original BML, we look at the asymptotic
speed of the system, i.e. the (asymptotic) average number of steps a
car moves in one turn. It is easily seen that this speed is the same
for all blue cars and all red cars, since the number of steps any
two cars of the same color have taken cannot differ by more than
$N$. It is somewhat surprising, perhaps, that these two speed must
also be the same. For instance, there is no configuration for which
a blue car completes 2 cycles for every 1 a red car completes.

\subsection{Time-normalized model} \label{normalization}

\paragraph{}
Though less natural, it will be sometimes useful to consider the
equivalent description, in which two rows of cars - red and blue -
are placed one beneath the other, with a "special place" - the
junction - where at most one car can be present at any time. In
every step first the red line shifts one to the right (except cars
immediately to the left of the junction, if it contains a blue car)
and then the blue line does the same. Furthermore, instead of having
the cars move to the right, we can have the junction move to the
left, and when a blue car is in the junction, the (possibly empty)
sequence of red cars immediately to the left of the junction moves
to the left, and vice verse. Figure \ref{timeimg} illustrates the
correspondence between these models.

\paragraph{}
From the discussion above we get the following equivalent system,
which we will call the time-normalized junction:

\begin{figure}
\centering
\includegraphics[scale=0.5]{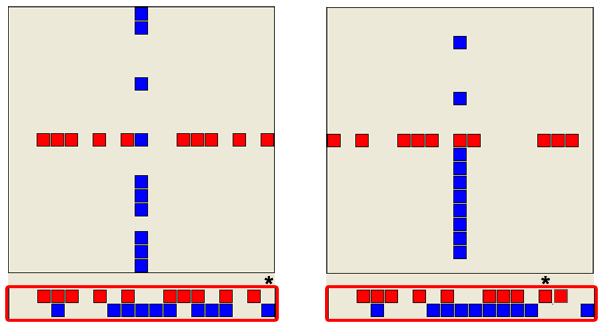}
\caption{On the left: a junction configuration and the analogous
configuration beneath it. The junction is marked with an asterisk.
On the right: same configuration after 3 turns in both views}
\label{timeimg}
\end{figure}

\begin{enumerate}
\item Fix $S=N-1 \in \Z_N$ and fix some initial configuration $\{ R_i \}_{i=0}^{N-1}, \{ B_i \}_{i=0}^{N-1} \in \{0,1\}^N$ representing the red and blue cars respectively,
i.e. $R_i=1$ iff there's a red car at the $i$-th place. We require
that $\sum R_i = \sum B_i = p \cdot N$, and at place $S$ ($=N-1$)
there is at most one car in both rows \footnote{The last requirement
is that the junction itself can contain only one car at the
beginning}.
\item In each turn:
\begin{itemize}
\item If place $S$ contains a blue car, and place $S-1$ contains a red car (if $B_S = R_{S-1} = 1$), push this car one step to the left. By pushing the car we mean also moving all red cars that immediately followed it one step to the left, i.e. set $R_{S-1} = 0$, $R_{S-i} = 1$ for $i = \min_{j \geq 1} [R_{S-j} = 0]$.
\item If place $S$ does not contain a blue car and place $S-1$ contains both a red and a blue car (if $B_S = 0$ and $R_{S-1} = B_{S-1} = 1$), push the blue car at $S-1$ one step to the left (set $B_{S-1} = 0$ and $B_{S-i} = 1$ for $i = \min_{j \geq 1} [B_{S-j} = 0]$).
\item set $S = S-1$
\end{itemize}
\end{enumerate}

Note that indeed the dynamics above guarantee that after a turn
place $(S-1)$ - the new junction - contains at most one car.
Generally, as long as cars flow freely, the time-normalized system
configuration does not change (except for the location of the
junction). Cars in the time-normalized system configuration actually
move only when some cars are waiting for the junction to clear (in
the non-time-normalized system).

\section{Analysis of an $(N,p)$ junction}
\paragraph{}
Our analysis of the junction will argue that for any $p$, regardless
of the initial configuration, the system will reach some optimal
speed, depending only on $p$. First let us state what is this
optimal speed:

\begin{theorem}
For a junction with density $p$, the maximal speed for any
configuration is $\min(1, \frac{1}{2p})$.

\label{maxspeedthm}

\end{theorem}

\paragraph{Proof}
Obviously, the speed cannot exceed 1. Observe the system when it
reaches its stable state, and denote the speed by $s$. Thus, at time
$t$, a car has advanced $ts(1+o(ts))$ steps (on average) which means
that it has passed the junction $ts/N (1+o(1))$ times. As only one
car can pass the junction at any time, the total number of cars
passing the junction at time $t$ is bounded by $t$. Therefore, we
have $2pN \times ts/N \le t$ which implies $s \le 1/2p$.
$\blacksquare$

\paragraph{}
We will now show that the system necessarily reaches these speeds
from any starting configuration.

\subsection{The case $p<0.5$}

\paragraph{}
We begin by proving that for $p < 0.5$ the junction will eventually
reach speed very close to 1. A system in a stable state is said to
be \emph{free-flowing} if no car is ever waiting to enter the
junction.

\begin{lemma}
A junction is free-flowing iff the time-normalized junction
satisfies:

(1) For all $0 \leq i \leq N-1$ there is only one car in place $i$ in both rows.

(2) For all $0 \leq i \leq N-1$ if place $i$ contains a blue car, place $(i-1) \mod N$

    does not contain a red car.

\label{violate1}
\end{lemma}

\paragraph{Proof}
Obviously, this is just a reformulation of free-flowing.
$\blacksquare$

\paragraph{}
We will now turn to show that for $p<0.5$, following the system
dynamics of the time-normalized junction will necessarily bring us
to a free flowing state, or at worst an "almost" free flowing state,
meaning a state for which the system speed will be arbitrarily close
to 1 for large enough $N$.

\paragraph{}
For this let us consider some configuration and look at the set of
"violations" to lemma \ref{violate1}, i.e. places that either
contain both a blue and a red car, or that contain a blue car and a
red car is in the place one to the left. As the following lemma will
show, the size of the set of violations is very closely related to
the system speed, and posses a very important property: it is
non-increasing in time.

\paragraph{}
More formally, For a configuration $R, B$ we define two disjoint
sets for the two types of violations: $$V_B = \{ 0 \leq i \leq N-1 :
R_{i-1} = B_{i-1} = 1, B_i = 0 \}$$
$$V_R = \{ 0 \leq i \leq N-1 : R_{i-1} = B_{i} = 1 \}$$
Also, let $V = V_B \cup V_R$ be the set of all violations. It will
be sometimes useful to refer to a set of indicators $X =
\{X(i)\}_{i=0}^{N-1}$, where each $X(i)$ is 1 if $i \in V$ and 0
otherwise, thus $|V| = \sum_{i=0}^{N-1} X(i)$.

For a junction with some initial configuration $R, B$, let $R^t,
B^t$ be the system configuration at time $t$, and let $V_t = V_{B^t}
\cup V_{R^t}$ be the set of violations for this configuration, and
$X^t$ be the corresponding indicator vector. Similarly, let $S^t$
denote the junction's position at time $t$.

\begin{lemma}
\
\begin{enumerate}
\item $|V_{t+1}| \leq |V_t|$
\item For any $t$, the system speed is at least $(1+\frac{|V_t|}{N})^{-1} \geq 1 - \frac{|V_t|}{N}$
\end{enumerate}
\label{violate2}
\end{lemma}

\paragraph{Proof}
Property (1) follows from the system dynamics. To see this, examine
the three possible cases for what happens at time $t$:

\begin{enumerate}
\item If in turn $t$ place $S^t$ does not contain a violation then the configurations do not change during the next turn, i.e. $R^{t+1} = R^t$, $B^{t+1}=B^t$ hence clearly $|V_{t+1}| = |V_t|$.

\item If $S^t \in V_{B^t}$ , then the configuration $B^{t+1}$ changes in two places:
\begin{enumerate}
\item $B^t_{S^t-1}$ is changed from 1 to 0. Thus, place $S^t$ is no longer in $V_{B^t}$, i.e. $X^t(S^t)$ changes from 1 to 0.
\item $B^t_{S^t-i}$ is changed from 0 to 1 for $i = \min_{j \geq 1} (B^t_{S^t-j} = 0)$. This may affect $X^{t+1}(S^t-i)$, and $X^{t+1}(S^t-i+1)$. However, for place $S^t-i+1$, by changing $B^t_{S^t-i}$ from 0 to 1 no new violation can be created (since by definition of $i$, $B^{t+1}_{S^t-i+1}=1$, so $X^{t+1}(S^t-i+1) = 1$ iff $R^{t+1}_{S^t-i} = 1$ regardless of $B^{t+1}_{S^t-i}$).
\end{enumerate}
For other indices $X^{t+1}(i) = X^t(i)$ since $R, B$ do not change, so between times $t$ and $t+1$, we have that $X^t(S^t)$ changes from 1 to 0, and at worst only one other place - $X^t(S^t-i+1)$ changes from 0 to 1, so $|V_{t+1}| = \sum_{i=0}^{N-1} X^{t+1}(i) \leq \sum_{i=0}^{N-1} X^t(i) = |V_t|$.

\item Similarly, if place $S^t \in V_{R^t}$ then the configuration $R^{t+1}$ changes in two places:
\begin{enumerate}
\item $R^{t+1}_{S^t-1}$ is changed from 1 to 0. Thus, $X^t(S^t)$ changes from 0 to 1.
\item $R_{S^t-i}$ is changed from 0 to 1 for $i = \min_{j \geq 1} (R^t_{S^t-j} = 0)$, affecting $X^{t+1}(S^t-i)$, $X^{t+1}(S^t-i+1)$. However for place $S^t-i$ changing $R_{S^t-i}$ does not affect whether this place is a violation or not, so at worst $X^t(S^t-i+1)$ changed from 0 to 1.
\end{enumerate}
By the same argument we get $|V_{t+1}| \leq |V_t|$.
\end{enumerate}
Therefore, $|V_{t+1}| \leq |V_t|$.
\paragraph{}
For property (2) we note that in the time-normalized system,
following a specific car in the system, its "current speed" is
$\frac{N}{N+k} = (1 + \frac{k}{N})^{-1}$ where $k$ is the number of
times the car was pushed to the left during the last $N$ system
turns. We note that if a car at place $j$ is pushed to the left at
some time $t$, by some violation at place $S^t$, this violation can
reappear only to the left of $j$, so it can push the car again only
after $S^t$ passes $j$. Hence any violation can push a car to the
left only once in a car's cycle (i.e. $N$ moves). Since by (1) at
any time from $t$ onwards the number of violations in the system is
at most $|V_t|$, then each car is pushed left only $|V_t|$ times in
$N$ turns, so its speed from time $t$ onwards is at least $(1 +
\frac{|V_t|}{N})^{-1}> 1 - \frac{|V_t|}{N}$ as asserted.
$\blacksquare$

\paragraph{}
With lemma \ref{violate2} at hand we are now ready to prove system
self organization for $p<0.5$. We will show that for $p<0.5$, after
$2N$ system turns $|V_t| = O(1)$, and hence deduce by part (2) of
lemma \ref{violate2} the system reaches speed $1 - O(\frac{1}{N})
\ra1$ (as $N\ra \infty$).

As the junction advances to the left it pushes some car sequences,
thus affecting the configuration to its left. The next lemma will
show that when $p<0.5$, for some $T<N$, the number of cars affected
to the left of the junction is only a constant, independent of $N$.

\begin{lemma}
Consider a junction with density $p<0.5$. There exists some constant
$C=C(p)=\frac{p}{1-2p}$, independent of $N$, for which:

From any configuration $R, B$ with junction at place $S$ there exist
some $0<T<N$ such that after $T$ turns:

(1) For $i\in \{S-T, \ldots, S\}$, $X^T(i)=0$ (i.e. there are no
violations there)

(2) For $i \in \{S+1, \ldots, N-1, 0, \ldots , S-T-C\}$, $R^T_i =
R^0_i$ and $B^T_i = B^0_i$ ($R, B$ are unchanged there) \label{ugly}
\end{lemma}

\paragraph{Proof}
First let us consider $T=1$. For $T=1$ to not satisfy the lemma
conclusions, there need to be a car sequence (either red or blue) of
length exceeding $C$, which is pushed left by the junction as it
moves. Progressing through the process, if, for some $T<N-C$, the
sequence currently pushed by the junction is shorter then $C$, then
this is the $T$ we seek. Therefore, for the conclusions \emph{not}
to hold, the length of the car sequence pushed by the junction must
exceed $C$ for all $0<T<N$. If this is the case, then leaving the
junction we see alternating red and blue sequences, all of lengths
exceeding $C$ and one vacant place after any blue sequence and
before any red one.

However, if this is the case for all $0<T\le T'$ then the average
number of cars per location in $\{S-T',\ldots,S\}$ at time $T'$ must
be at least $\frac{2C}{2C+1}$ (at least $2C$ cars between vacant
places). Therefore, the total number of cars in $\{S-T',\ldots,S\}$
at time $T'$ is more than $\frac{2C}{2C+1} T'=2p T'$ (Recall that
$C=\frac{2p}{1-2p}$).

Since there are only $2pN$ cars in the system, this cannot hold for
all $T$ up to $N$. Thus, there must be some time for which the
conclusions of the lemma are satisfied.$\blacksquare$

\paragraph{}
We are now ready to easily prove the main result for this section.

\begin{theorem}
A junction of size $N$ with density $p<0.5$ reaches speed of
$1-\frac{C(p)}{N}$ from any initial configuration. \label{speed1thm}
\end{theorem}

\paragraph{Proof}
Let $R, B$ be some initial configuration, with $S=N-1$ and let $V$
the corresponding set of violations and $X$ the matching indicators
vector. By lemma \ref{ugly} there exist $T_0>0$ for which
$X^{T_0}(i) = 0$ for $i \in [N-1-T_0, N-1]$. Now starting at
$R^{T_0}, B^{T_0}$ and $S=N-1-T_0$ reusing the lemma there exist
$T_1>0$ s.t. $X^{T_0 + T_1}(i) = 0$ for $i \in [N-1-T_1, N-1-T_0]$,
and also, as long as $N-1-T_0-T_1 > C(p) = \frac{2p}{1-2p}$, $X^{T_0
+ T_1}(i) = X^{T_0}(i) = 0$ for $i \in [N-1-T_0, N-1]$ as well.

Proceeding in this manner until $T = \sum T_i \geq N$ we will get
that after $T$ turns, $X^T(i) = 0$ for all but at most $C(p)$
places, hence by lemma \ref{violate2} the system speed from this
time onward is at least $1-\frac{C(p)}{N}$. $\blacksquare$

\begin{figure}
\centering
\includegraphics[scale=0.5]{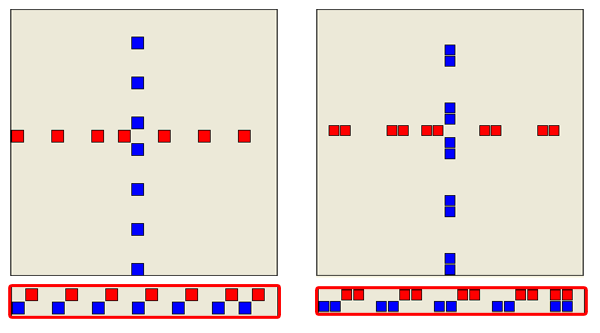}
\caption{(a) A junction configuration with density $p=\frac{1}{3}$
($\frac{p}{1-2p} = 1$) not reaching speed 1 (speed is $1 -
\frac{1}{N}$). (b) A similar construction for $p=0.4$
($\frac{p}{1-2p} = 2$) reaching speed of only $1 - \frac{2}{N}$}
\label{counterimg}
\end{figure}

\paragraph{}
We remark one cannot prove an actual speed 1 (rather than
$1-o(1)$) for the case $p<\frac{1}{2}$ since this is definitely
not the case. Figure \ref{counterimg} (a) demonstrates a general
construction of a junction with speed $1-\frac{1}{N}$ for all $N$.
Finally we remark that a sharper constant $C'(p)$ can be provided,
that will also meet a provable lower bound for the speed (by
constructions similar to the above). This $C'(p)$ is exactly half
of $C(p)$. We will see proof of this later in this paper as we
obtain theorem \ref{speed1thm} as a corollary of a different
theorem, using a different approach.

\subsection{Number of segments for $p<0.5$}

\paragraph{}
The proof in the previous section were combinatorial in nature and
showed that for a density $p<0.5$ the system can slow-down only by
some constant number independent of $N$ regardless of the
configuration. As we turn to examine the properties of the stable
configuration of course we can not say much if we allow any initial
configuration. There are configurations in which the number of
different segments in each row will be $\Theta(N)$ while clearly if
the cars are arranged in a single red sequence and a single blue
sequence in the first place, we will have only one sequence of each
color at any time.

\paragraph{}
However, we will show that for a random initial configuration, the
system will have linearly many different segments of cars with high
probability.

\begin{theorem}
A junction of size $N$ with density $p<0.5$, started from a random
initial configuration, will have $\Theta(N)$ different segments at
all times with high probability (w.r.t. $N$). \label{segments1thm}
\end{theorem}

\paragraph{Proof}
As we already seen in the proof of \ref{ugly}, as the system
completes a full round, since every place (but a constant number
of places) contains at most a single car there must be $(1-2p)N$
places in which no car is present. Each two such places that are
not adjacent must correspond to a segment in the cars
configuration.

\paragraph{}
It is evident by the system dynamics, that the number of places
for which $R_i = B_i = R_{i-1} = B_{i-1} = 0$ is non-increasing.
More precisely, only places for which $R_i = B_i = R_{i-1} =
B_{i-1} = 0$ in the initial configuration can satisfy this in the
future. In a random initial configuration with density $p$, the
initial number of these places is expected to be $(1-p)^4N$, and
by standard CLT, we get that with high probability this number is
at most $((1-p)^4 + \varepsilon)N$. Thus, the number of different
segments in the system configuration at any time is at least
$((1-2p) - (1-p)^4 - \varepsilon)N$. However for $p$ very close to
$\frac{1}{2}$ this bound may be negative, so this does not
suffice.

\paragraph{}
To solve this we note that similarly also for any fixed $K$, the
number of consecutive $K$ empty places in a configuration is
non-increasing by the system dynamics, and w.h.p. is at most
$((1-p)^{2K} + \varepsilon)N$ for an initial random state. But
this guarantees at least $\frac{(1-2p) - (1-p)^{2K} -
\varepsilon}{K-1}N$ different segments in the system. Choosing
$K(p)$ (independent of N) sufficiently large s.t. $(1-p)^{2K} <
(1-2p)$ we get a linear lower bound on the number of segments from
a random initial configuration. $\blacksquare$

\subsection{The case $p>0.5$}

\paragraph{}
The proofs of speed optimality and segment structure for $p > 0.5$
will rely mainly of a the combinatorial properties of a \emph{stable
configuration}. A \emph{stable configuration} for the system is a
configuration that re-appears after running the system for some $M$
turns. Since for a fixed $N$ the number of possible configurations
of the system is finite, and the state-transitions (traffic rules)
are time independent, the system will necessarily reach a stable
configuration at some time regardless of the starting configuration.

\paragraph{}
We will use mainly two simple facts that must hold after a system
reached a stable configurations: (a) $|V_t|$ cannot change - i.e. no
violation can disappear. This is clear from lemma \ref{violate2};
(b) Two disjoint segments of car cannot merge to one (i.e. one
pushed until it meets the other), since clearly the number of
segments in the system is also non-increasing in time.

\paragraph{}
These two facts alone already provide plenty of information about
the stable configuration for $p>0.5$. We begin with the following
twin lemmas on the stable state.

\begin{figure}
\centering
\includegraphics[scale=0.55]{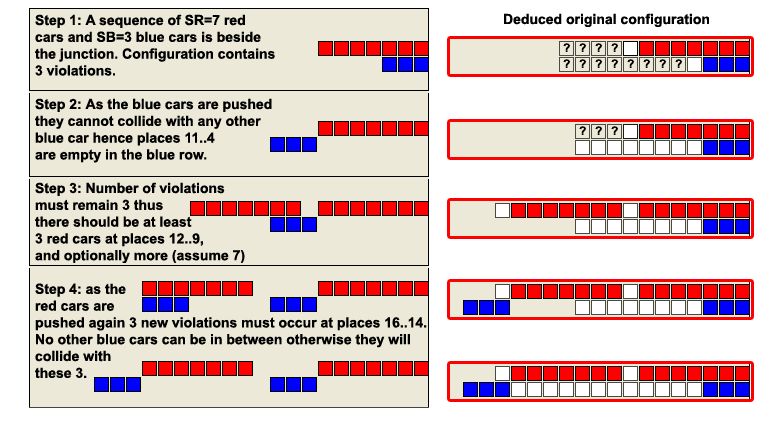}
\caption{Sketch of proof ideas for lemmas \ref{uglier1} (steps 1-3) and \ref{uglier2} (step 4)}
\label{proofimg}
\end{figure}

\begin{lemma}
Let $R,B$ be a stable configuration with junction at $S=0$ and
$B_0=0$. Assume that there is a sequence of exactly $s_R$
consecutive red cars at places $[N-s_R, N-1]$ and $s_B$ blue cars at
places $[N-s_B, N-1]$, $s_R, s_B \geq 1$. Then:
\begin{enumerate}
\item $B_i = 0$ for $i \in [N-s_R-s_B-1, N-s_B-1]$
\item $R_i = 1$ for $i \in [N-s_R-s_B-1, N-\max(s_R,s_B)-2]$
\item $R_i = 0$ for $i \in [N-\max(s_R,s_B)-1, N-s_R-1]$
\end{enumerate}
\label{uglier1}
\end{lemma}
\paragraph{Proof}
Since it is easy to lose the idea in all the notations, a visual
sketch of the proof is provided in figure \ref{proofimg}.

\paragraph{}
By the assumptions:

- $B_i = 1$ for $i \in [N-s_B, N-1]$, $B_{N-s_B-1} = 0$

- $R_i = 1$ for $i \in [N-s_R, N-1]$, $R_{N-s_R-1} = 0$

\paragraph{}
To get (1), we note that with the assumption $B_0 = 0$, the blue
sequence will be pushed to the left in the next $s_R$ turns, so by
restriction (b), $B_i = 0$ for $i \in [N-s_B-s_R-1, N-s_B-1]$ since
otherwise 2 disjoint blue segments will merge while the sequence is
pushed.

\paragraph{}
Thus following the system, after $s_R$ turns we will get: $B_i = 1$
for $i \in [N-s_R-s_B,N-s_R-1]$ and $B_i = 0$ for $i \in
[N-s_R,N-1]$, and $B_i$ not changed left to $N-s_R-s_B$, and $R$
unchanged.
\paragraph{}
Note that originally $R,B$ contained $\min(s_R,s_B)$ consecutive
violations in places $[N-\min(s_R,s_B), 0]$ which all vanished after
$s_R$ turns. Possible violations at places $[N-\max(s_R,s_B),
N-\min(s_R,s_B)]$ remained as they were. From here we get that we
must have $R_i=1$ for $\min(s_R,s_B)$ places within $[N-s_R-s_B-1,
N-\max(s_R,s_B)-1]$. Since for place $N-\max(s_R,s_B)-1$ either $R$
or $B$ are empty by the assumption, we must therefore have $R_i = 1$
for $i \in [N-s_B-s_R-1, N-\max(s_R,s_B)-2]$, giving (2).
\paragraph{}
If we follow the system for $s_B$ more steps we note that any red
car in $[N-\max(s_R,s_B)-1, N-s_R-1]$ will be pushed left until
eventually hitting the red car already proven to be present at
$N-\max(s_R,s_B)-2]$, thus $R_i = 0$ for $i \in [N-\max(s_R,s_B)-1,
N-s_R-1]$ giving (3). $\blacksquare$

\paragraph{}
The following lemma is completely analogous when reversing the
roles of $R,B$, and can be proven the same way.

\begin{lemma}
Let $R,B$ be a stable configuration with junction at $S=0$ and
$B_0=1$. Assume that there is a sequence of exactly $s_R$
consecutive red cars at places $[N-s_R, N-1]$ and $s_B$ blue cars at
places $[N-s_B+1, 0]$, $s_R, s_B \geq 1$. Then:
\begin{enumerate}
\item $R_i = 0$ for $i \in [N-s_B-s_R-1, N-s_R-1]$
\item $B_i = 1$ for $i \in [N-s_B-s_R, N-\max(s_B,s_R)-1]$
\item $B_i = 0$ for $i \in [N-\max(s_B,s_R),N-s_B]$
\end{enumerate}
\label{uglier2}
\end{lemma}

\paragraph{}
Putting lemmas \ref{uglier1}, \ref{uglier2} together we get the
following characterization for stable configurations:
\begin{lemma}
Let $R,B$ be a stable configuration with junction at $S=0$ and
$B_0=0$. Assume that there is a sequence of exactly $s_R$
consecutive red cars at places $[N-s_R, N-1]$ and $s_B$ blue cars at
places $[N-s_B, N-1]$. Denote $M=\max(s_B, s_R)$ Then
\begin{enumerate}
\item There are no additional cars are at $[N-M, N-1]$.
\item Place $i = N-M-1$ is empty - i.e. $R_i = B_i = 0$
\item Starting at $N-M-2$ there is a sequence of $K_1 \geq \min(s_B, s_R)$ places for which $R_i = 1; B_i = 0; (i \in [N-M-K_1-1 , N-M-2])$
\item Starting $N-M-K_1-2$ (i.e. right after the red sequence) there is a sequence of $K_2 \geq \min(s_B, s_R)$ places for which $B_i = 1; R_i = 0; (i \in [N-M-K_1-K_2-1, N-M-K_1-2])$.
\end{enumerate}
\label{ugliest}
\end{lemma}

\begin{figure}
\includegraphics[scale=0.6]{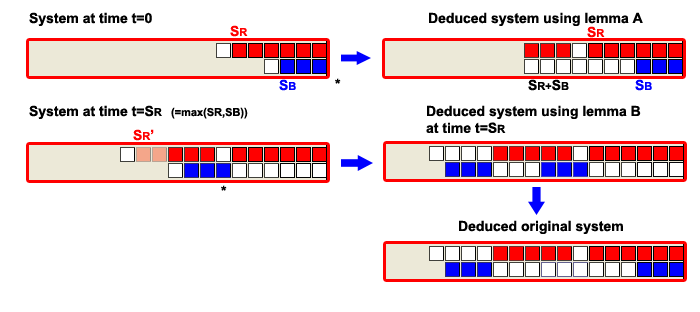}
\caption{Lemmas \ref{uglier1} (A) and \ref{uglier2} (B) combined together yield lemma \ref{ugliest}}
\label{proof3img}
\end{figure}

\paragraph{Proof}
Figure \ref{proof3img} outlines the proof, without the risk of
getting lost in the indices. For the full proof, first get from
lemma \ref{uglier1} that:
\begin{enumerate}
\item $B_i = 0$ for $i \in [N-s_R-s_B-1, N-s_B-1] $
\item $R_i = 1$ for $i \in [N-s_R-s_B-1, N-\max(s_R,s_B)-2]$
\item $R_i = 0$ for $i \in [N-\max(s_R,s_B)-1, N-s_R-1]$
\end{enumerate}

\paragraph{}
From (1),(3) we get in particular, $B_i = 0$ for $i \in
[N-\max(s_R,s_B)-1, N-s_B-1]$, and $R_i = 0$ for $i \in
[N-\max(s_R,s_B)-1, N-s_R-1]$, so indeed no additional cars are at
$[N-\max(s_B,s_R), N-1]$, and place $N-\max(s_R,s_B)-1$ is empty,
proving claims 1,2 in the lemma.

\paragraph{}
From (2),(3) we get $B_i = 0$ for $i \in [N-s_R-s_B-1,
N-\max(s_R,s_B)-2]$ and $R_i = 1$ for $i \in [N-s_R-s_B-1,
N-\max(s_R,s_B)-2]$, thus places $[N-s_R-s_B-1, N-\max(s_R,s_B)-2]$
contain a sequence of length $\min(s_R,s_B)$ of red cars with no
blue cars in parallel to it. This sequence is possibly a part of a
larger sequence of length $s_R' \geq \min(s_R,s_B)$, located at
$[N-\max(s_R,s_B)-s_R'-1, N-\max(s_R,s_B)-2]$.

\paragraph{}
Now running the system for $\max(s_R,s_B)$ turns, we will have the
junction at place $S=N-\max(s_R,s_B)-1$, $B_S = 1$, followed by
sequences of $s_R'$ reds and $s_B' = \min(s_B, s_R) (\leq s_R')$
blues. Applying lemma \ref{uglier2} for the system (rotated by
$\max(s_R,s_B)$ i.e. for $N' = N-\max(s_R,s_B)-1$):
\begin{enumerate}
\item $R_i = 0$ for $i \in [N'-s_B'-s_R'-1, N'-s_R'-1] = [N-\max(s_R,s_B)-\min(s_R, s_B)-s_R'-2, N-\max(s_R,s_B)-s_R'-2] = [N-s_R-s_B-s_R'-2, N-\max(s_R,s_B)-s_R'-2]$
\item $B_i = 1$ for $i \in [N'-s_B'-s_R', N'-\max(s_B',s_R')-1] = [N-s_R-s_B-s_R'-1, N-\max(s_R,s_B)-\max(s_B',s_R')-2] = [N-s_R-s_B-s_R'-1, N-\max(s_R,s_B)-s_R'-2]$
\item $B_i = 0$ for $i \in [N'-\max(s_B',s_R'),N'-s_B'] = [N-\max(s_R,s_B)-s_R'-1, N-\max(s_R,s_B)-\min(s_R,s_B)-1] = [N-\max(s_R,s_B)-s_R'-1, N - s_R - s_B - 1]$
\end{enumerate}

\paragraph{}
In particular from (3) we get that no blue cars are in parallel to
the entire red segment in $[N-\max(s_R,s_B)-s_R'-1,
N-\max(s_R,s_B)-2]$: We were previously assured this is true up to
place $N-s_R-s_B-1$, and for places $[N-\max(s_R,s_B)-s_R'-1,
N-s_R-s_B-2] \subseteq [N-\max(s_R,s_B)-s_R'-1,
N-\max(s_R,s_B)-s_B'-1]$ this holds by (3).

\paragraph{}
Furthermore by (2) we get that a sequence of blue cars which begins
from place $N-\max(s_R,s_B)-s_R'-2$ with no red cars in parallel to
it by (1). Note that $N-\max(s_R,s_B)-s_R'-2$ is exactly to the left
of  $N-\max(s_R,s_B)-s_R'-1$ where the red sequence ended. Now
clearly choosing $K_1 = s_R'; K_2 = \min(s_R, s_B)$ we get claims
3,4,5 in the lemma. $\blacksquare$

\paragraph{}
Putting it all together we can now get a very good description of a stable state.

\begin{theorem}
Let $R,B$ be a stable configuration with junction at $S=0$ and
$B_0=0$. Assume that there is a sequence of exactly $s_R$
consecutive red cars at places $[N-s_R, N-1]$ and $s_B$ blue cars at
places $[N-s_B, N-1]$. Denote $M=\max(s_B,s_R)$.

Then no additional cars are at $[N-M, N-1]$, and at places $[0,N-M-1]$ the configurations $R,B$ satisfies:
\begin{enumerate}
\item Each place contains at most one type of car, red or blue.
\item Place $N-M-1$ is empty. Each empty place, is followed by a sequence of places containing red cars immediately left to it, which is followed by a sequence of places containing blue cars immediately left to it.
\item Any sequence of red or blue cars is of length at least $\min(s_R, s_B)$
\end{enumerate}
\label{stablethm}
\end{theorem}

\begin{figure}
\includegraphics[scale=0.6]{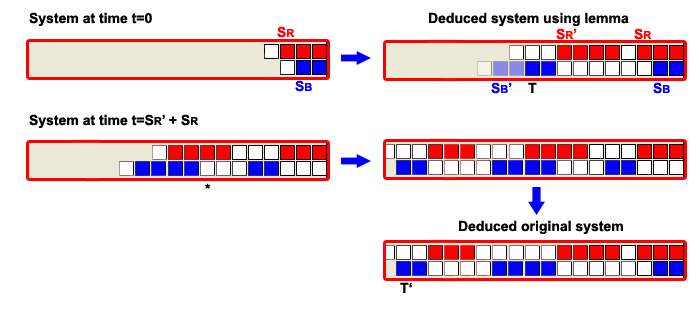}
\caption{Repetitively applying lemma \ref{ugliest}, we unveil a
longer segment in the configuration, $[T,N-1]$, for which
properties of theorem \ref{stablethm} hold} \label{proof4img}
\end{figure}

\paragraph{Proof}
This is merely applying lemma \ref{ugliest} repeatedly. Applying
lemma \ref{ugliest} we know that there exist $K_1, K_2 \geq
\min(s_R, s_B)$ such that: Place $N-M-1$ is empty, followed by $K_1$
consecutive places with only red cars and $K_2$ consecutive places
with only blue cars left to it, thus the assertion holds for the
segment $[T, N-M-1]$ for $T=N-M-1-K_1$.
\paragraph{}
Now we completely know $R,B$ in $[N-M-1-K_1,N-M-1]$, and this is
enough to advance the system for $K_1 + M$ turns. The $s_B$ blue
segment is pushed left $s_R$ places, further pushing the $K_1$ red
sequence $\min(s_R, s_B)$ places to the left such that its last
$\min(s_R, s_B)$ cars now overlap with the $K_2$ blue sequence.
\paragraph{}
So after $K_1 + M$ turns the system evolves to a state where $S =
N-K_1-M$, $B_S=0$, and left to $S$ there are $K_2' > K_2$
consecutive blue cars and exactly $\min(s_R, s_B)$ consecutive red
cars. Noting that this time $M'=max(K_2', \min(s_R,s_B))=K_2'$, once
again we can deduce from lemma \ref{ugliest} that: there are no
additional cars in $[N-M-K_1-K_2',N-M-K_1-1]$ (thus we are assured
that the entire blue segment of length $K_2'$ does not have red cars
parallel to it), Place $N-M-K_1-K_2'-1$ is empty, followed by some
$K_3$ consecutive places with only red cars and $K_4$ consecutive
places with only blue cars left to it, for $K_3, K_4 \geq
\min(K_1,K_2) = \min(s_R,s_B)$ thus the assertion holds for the
segment $[T', N-M-1]$ for $T'=N-M-K_1-K_2'-K_3-1 < T$.
\paragraph{}
Repeatedly applying lemma \ref{ugliest} as long as $T>0$, we
repeatedly get that the assertion holds for some $[T,N-M-1]$ for $T$
strictly decreasing, so the assertion holds in $[0,N-M-1]$
$\blacksquare$

\begin{figure}
\centering
\includegraphics{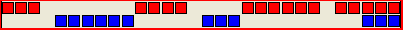}
\caption{A typical stable configuration}
\label{typical}
\end{figure}

\paragraph{}
We have worked hard for theorem \ref{stablethm}, but it will soon turn to be worthwhile.
Let us obtain some useful corollaries.

\begin{corollary}
Let $R,B$ be a stable configuration with junction at $S=0$ and
$B_0=0$. Assume that there is a sequence of exactly $s_R$
consecutive red cars at places $[N-s_R, N-1]$ and $s_B$ blue cars at
places $[N-s_B, N-1]$. Denote $m = \min(s_B,s_R)$ Then:
\begin{enumerate}
\item The number of blue segments in the system equals the number of red segments.
\item System speed is at least $(1+\frac{m}{N})^{-1}$
\item Total number of cars in the system is at most $N+m$
\item Total number of cars in the system is at least $\frac{2m}{2m+1}N+m$
\item Total number of segments in the system is at most $\frac{N}{m}+1$
\end{enumerate}
\label{stablecor}
\end{corollary}

\paragraph{Proof}
(1) Follows by the structure described in \ref{stablethm}, since
there is exactly one red and one blue sequence in $[N-M,N-1]$ (as
always $M = \max(s_R,s_B)$) and an equal number of reds and blues in
$[0, N-M-1]$ since any red sequence is immediately proceeded by a
blue sequence.

\paragraph{}
For (2) we note that by \ref{stablethm} the configuration $R,B$
contains exactly $m$ violations : the $m$ overlapping places of the
segments $s_R, s_B$ in $[N-M,N-1]$ are the only violations in the
configuration since any places in $[0,N-M-1]$ contains a single car,
and no red car can be immediately to the left of a blue car. By
\ref{violate2} the system speed is hence at least
$(1+\frac{m}{N})^{-1}$.

\paragraph{}
By \ref{stablethm} any place in $[0,N-m-1]$ contains at most one
car, or it is empty, and places $[N-m,N-1]$ contain both a red car
and a blue car. Thus the total number of cars in the system is at
most $N-m+2m =N+m$ giving (3).

\paragraph{}
On the other hand, Since any sequence of a red car or a blue car is
of length at least $m$, an empty place in $[0, N-m]$ can occur only
once in $2m+1$ places, and other places contain one car. Thus the
number of cars is lower bounded by $\frac{2m}{2m+1}(N-m) + 2m =
\frac{2m}{2m+1}N + \frac{2m^2+2m}{2m+1} \geq \frac{2m}{2m+1}N + m$
giving (4).

\paragraph{}
The last property follows the fact that any sequence of cars is of
length at least $m$, and by (3) total number of cars is at most
$N+m$, thus number of different sequences is at most $\frac{N+m}{m}
= \frac{N}{m}+1$. $\blacksquare$

\begin{theorem}
All cars in the system have the same asymptotic speed.
\end{theorem}

\paragraph{Proof}
As we have seen, when it reaches a stable state, the system consists
of alternating red and blue sequences of cars. Obviously, the order
of the sequences cannot change. Therefore, the difference between
the number of steps two different cars have taken cannot be more
then the length of the longest sequence, which is less then $N$.
Thus, the asymptotic (w.r.t. $t$) speed is the same for all cars.
$\blacksquare$

\paragraph{}
With these corollaries we can now completely characterize the stable
state of a junction with $p > 0.5$, just by adding the final simple
observation, that since the number of cars in the model is greater
than $N$, there are violations at all times, including after
reaching a stable state. Now let us look at some time when the
junction reaches a violation when the system is in stable state. At
this point the conditions of theorem \ref{stablethm} are satisfied,
thus:

\begin{theorem}
A junction of size $N$ and density $p>0.5$ reaches speed of
$\frac{1}{2p} - O(\frac{1}{N})$ (i.e. arbitrarily close to the
optimal speed of $\frac{1}{2p}$, for large enough $N$), and contains
at most a bounded number (depending only on $p$) of car sequences.
\label{largepthm}
\end{theorem}

\paragraph{Proof}
We look at the system after it reached a stable state. Since $2pN >
N$ at some time after that conditions of theorem \ref{stablethm} are
satisfied for some $s_R, s_B \geq 1$. Let $m = \min(s_R,s_B)$ at
this time. Using claims (3),(4) in corollary \ref{stablecor} we get:
$$\frac{2m}{2m+1}N + m \leq 2pN \leq N + m$$
From here we get
$$(*)\ \ (2p-1)N \leq m \leq (2p - \frac{2m}{2m+1})N = (2p - 1)N + \frac{N}{2m+1} $$
and reusing $m \geq (2p-1)N$ on the left hand size we get:
$$(2p-1)N \leq m \leq (2p-1)N + \frac{1}{4p-1}$$
For $C = \frac{1}{4p-1}$ a constant independent of $N$ ($C=O(1)$).
So $m = (2p-1)N + K$, for $K \leq C$. Now by claim (2) in
\ref{stablecor}, system speed is at least
$$(1+\frac{m}{N})^{-1} = (1+\frac{(2p-1)N + K}{N})^{-1} = (2p + \frac{K}{N})^{-1} \geq_{(2p > 1)} \frac{1}{2p} - \frac{K}{N}$$
But by theorem \ref{maxspeedthm} system speed is at most
$\frac{1}{2p}$, thus system speed is exactly $\frac{1}{2p} -
\frac{K'}{N}$ for some $0 \leq K' \leq K \leq C = O(1)$, thus $K' =
O(1)$ proving the first part of the theorem.

\paragraph{}
By (5) we get total number of segments in the system is at most
$\frac{N}{m}+1$, applying $m \geq (2p-1)N$ we get the number of
segments is bounded by:
$$\frac{N}{m}+1 \leq \frac{1}{2p-1}+1 = \frac{2p}{2p-1} = O(1)$$
Thus the second part proven. $\blacksquare$

\subsection{$p < 0.5$ revisited}
\paragraph{}
The characterization in $\ref{stablethm}$ can be also proven useful
to handle $p < 0.5$. Actually the main result for $p<0.5$, theorem
$\ref{speed1thm}$, can be shown using similar technique, and even
sharpened.\footnote{To theorem \ref{speed1thm} defence, we should
note that key components of its proof, such as lemma \ref{violate2},
were also used in the proof of $\ref{stablethm}$}

\begin{corollary}
For $p < 0.5$ the junction reaches speed of at least $1 -
\frac{C(p)}{N}$, for $C(p) = \left\lfloor
\frac{p}{1-2p}\right\rfloor$. In particular, for $p < \frac{1}{3}$
the junction reaches speed 1, for any initial configuration.
\label{speed1cor}
\end{corollary}

\paragraph{Proof}
Let $R,B$ be any initial configuration. Looking at the configuration
after it reached the stable state, if the system reached speed 1 we
have nothing to prove. Assume the speed is less than 1. Since in
this case violations still occur, at some time the stable
configuration will satisfy theorem \ref{stablethm}. As before,
letting $m = \min(s_R,s_B)$ at this time, by claim (4) in corollary
\ref{stablecor} we have:
$$\frac{2m}{2m+1}N + m \leq 2pN  \ \Rightarrow m \leq (2p-1+\frac{1}{2m+1})N$$
In particular $2p-1+\frac{1}{2m+1} > 0$, rearranging we get $m < \frac{p}{1-2p}$, and since $m \in \Z$, $m \leq \left\lfloor \frac{p}{1-2p}\right\rfloor = C(p)$. $m$ must be positive, thus for $p<\frac{1}{3}$, having $C(p) = \left\lfloor \frac{p}{1-2p}\right\rfloor = 0$ we get a contradiction, thus the assumption (speed $<1$) cannot hold. For $p\geq \frac{1}{3}$ by claim (2) in corollary \ref{stablecor} the system speed is at least $(1+\frac{C(p)}{N})^{-1} \geq 1 - \frac{C(p)}{N}$. $\blacksquare$

\subsection{The critical $p=0.5$}
\paragraph{}
Gathering the results so far we get an almost-complete description
of the behaviour of a junction. For junction of size $N$ and density
$p$:
\begin{itemize}
\item If $p < 0.5$ the junction will reach speed $1 - o(1)$ (asymptotically optimal), and contain linearly many different segments in the stable state.
\item If $p > 0.5$ the junction will reach speed $\frac{1}{2p} - o(1)$ (asymptotically optimal), and contain constant many segments in the stable state.
\end{itemize}

\paragraph{}
From the description above one sees that the junction system goes
through a sharp phase transition at $p=0.5$, as the number of
segments of cars as the system stabilizes drops from being linear
to merely constant. The last curiosity, is what happens at
$p=0.5$. Once again by using the powerful theorem
$\ref{stablethm}$ we can deduce:

\begin{theorem}
A junction of size $N$ with $p=0.5$ reaches speed of at least $1 -
\frac{1}{\sqrt{N}}$ and contains at most $\sqrt{N}$ different
segments.
\end{theorem}

\paragraph{Proof}
For $p=0.5$ we have exactly $N$ cars in the system. As we reach
stable state, violations must still occur (since a system with
exactly $N$ cars must contain at least one violation), thus at some
time theorem \ref{stablethm} is satisfied. For $m = \min(s_R,s_B)$
at this time, by claim (4) in corollary \ref{stablecor} we have:
$$\frac{2m}{2m+1}N + m \leq N \ \Rightarrow m(2m+1) \leq N$$
Thus $m < \sqrt{N}$. From here by claim (2) the system speed is at
least $(1+\frac{\sqrt{N}}{N})^{-1} \geq 1 - \frac{1}{\sqrt{N}}$.
\paragraph{}
If $S$ is the number of segments, then by theorem \ref{stablethm}
we can deduce that the total number of cars, $N$, is: $2m$ cars in
places $[N-m, N-1]$, and $N-m-S$ cars in places $[0,N-M-1]$ (since
each place contains one car, except transitions between segments
that are empty).
$$N = (N-m-S) + 2m = N+m-S\ \Rightarrow S = m \leq \sqrt{N}$$
Thus the configuration contains at most $\sqrt{N}$ segments. $\blacksquare$

\paragraph{}
Simulation results show that these bounds are not only tight, but
typical, meaning that a junction with a random initial configuration
with density $p$ indeed has $O(\sqrt{N})$ segments in the stable
state, with the largest segment of size near $n^{1/2}$. This
suggests that the system undergoes a second order phase transition.
\newpage
\section{Simulation results}
\paragraph{}
Following are computer simulation results for the junction for
critical and near critical $p$, demonstrating the phase transition.
The columns below consist of $N$, the average asymptotic speed, the
average number of car segments in the stable state, the average
longest segment in the stable state, and the average number of
segments divided by $N$ and $\sqrt{N}$.

\subsection{$p = 0.48$}
\paragraph{}
For large $N$, system reaches speed $1$ and the average number of
segments is linear (approx. $0.037N$).
\paragraph{}

\begin{tabular}{|c|c|c|c|c|c|}
\hline
$N$ & $Speed$ & $No.\ segs$ & $Longest$ & $No.\ segs/N$ & $No.\ segs/\sqrt{N}$ \\
\hline
   1000 & 0.99970 &    38.7 &     6.8 & 0.0387 & 1.2238 \\
   5000 & 1.00000 &   186.4 &     8.5 & 0.0373 & 2.6361 \\
  10000 & 1.00000 &   369.8 &     7.5 & 0.0370 & 3.6980 \\
  50000 & 1.00000 &  1850.6 &     8.2 & 0.0370 & 8.2761 \\
\hline
\end{tabular}

\subsection{$p = 0.52$}
\paragraph{}
For large $N$, system reaches speed $0.961 = \frac{1}{2p}$ and the
average number of segments is about constant.
\paragraph{}
\begin{tabular}{|c|c|c|c|c|c|}
\hline
$N$ & $Speed$ & $No.\ segs$ & $Longest$ & $No.\ segs/N$ & $No.\ segs/\sqrt{N}$ \\
\hline
   1000 & 0.95703 &     5.7 &    76.7 & 0.0057 & 0.1802 \\
   5000 & 0.96041 &     6.9 &   330.0 & 0.0014 & 0.0976 \\
  10000 & 0.96091 &     7.3 &   416.1 & 0.0007 & 0.0730 \\
  50000 & 0.96142 &     7.2 &  3207.1 & 0.0001 & 0.0322 \\
 \hline
\end{tabular}

\subsection{$p = 0.5$}
\paragraph{}
At criticality, the speed is approaching 1 like $1 -
\frac{C}{\sqrt{N}}$ and the average number of segments is around
$0.43\cdot \sqrt{N}$.
\paragraph{}
\begin{tabular}{|c|c|c|c|c|c|}
\hline
$N$ & $Speed$ & $No.\ segs$ & $Longest$ & $No.\ segs/N$ & $No.\ segs/\sqrt{N}$ \\
\hline
   1000 & 0.98741 &    13.4 &    38.4 & 0.0134 & 0.4237 \\
   5000 & 0.99414 &    30.0 &    82.8 & 0.0060 & 0.4243 \\
  10000 & 0.99570 &    43.8 &   142.1 & 0.0044 & 0.4375 \\
  50000 & 0.99812 &    95.0 &   248.4 & 0.0019 & 0.4251 \\
\hline
\end{tabular}
\newpage
\section{Summary}
\paragraph{}
The fascinating phenomena observed in the BML traffic model are
still far from being completely understood. In this paper we showed
a very simplified version of this model, which, despite its relative
simplicity, displayed very similar phenomena of phase transition at
some critical density and of self-organization, which in our case
both can be proven and well understood.
\paragraph{}
We used two approaches in this paper: The first one was to use some
sort of "independence" in the way the system evolves, as in the way
we handle $p < 0.5$. We showed that the system self-organizes
"locally" and with a bounded affect on the rest of the
configuration, thus it will eventually organize globally. The second
approach is the notion of the stable configuration, i.e. we
characterize the combinatorial structure of any state that the
system can "preserve" and use it to show it is optimal (in a way
saying, that as long as we are not optimal the system must continue
to evolve).
\begin{figure}
\centering
\includegraphics[scale=0.6]{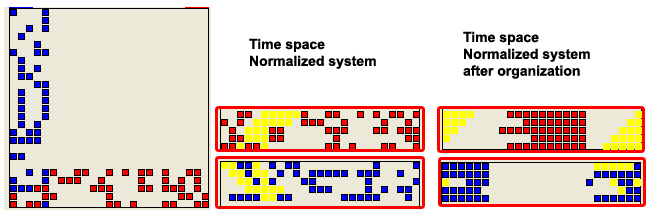}
\caption{K-Lanes junction in a time-space normalized view, before and after organization. Junction is marked in yellow and advances left each turn}
\label{klanesimg}
\end{figure}
\paragraph{}
Can these results be extended to handle more complicated settings
than the junction? Possibly yes. For example, considering a
$k$-lanes junction (i.e. $k$ consecutive red lines meeting $k$
consecutive blue rows), one can look at a time-space-normalized
version of the system, as shown in figure \ref{klanesimg}, with the
junction now being a $k\times k$ parallelogram traveling along the
red and blue lines, and "propagating violations" (which now have a
more complicated structure depending on the $k\times k$
configuration within the junction). Stable states of this
configurations seem to have the same characteristics of a single
junction, with a red (blue) car equivalent to some red (blue) car in
any of the $k$ red (blue) lines. Thus a zero-effort corollary for a
$k$-lanes junction is that it reaches speed 1 for $p <
\frac{1}{2k}$, but for $k$ nearing $O(N)$ this bound is clearly
non-significant. It is not surprising though, since combinatorics
alone cannot bring us too far, at least for the complete BML model
-- even as little as $2N$ cars in an $N^2$ size torus can be put in
a stuck configuration - i.e. reach speed 0.

\end{document}